\documentclass{article}
\usepackage[utf8]{inputenc}
\usepackage[
backend=biber,
style=alphabetic,
]{biblatex}
\usepackage{graphicx}
\usepackage{amsfonts}
\usepackage{amsmath}
\usepackage{amssymb}
\usepackage{amsthm}
\usepackage{float}
\newtheorem{theorem}{Theorem}[section]

\newtheorem{lemma}[theorem]{Lemma}

\title{Nonexistence of Perfect Euler Boxes with Semiprime Sides}
\author{Riley Tao}
\date{October 2021}

\begin{document}

\maketitle

\begin{abstract}
    In this paper we prove that there cannot exist a perfect Euler box with a semiprime side. We first display the proof, which uses nothing more than elementary number theory. Due to the elementary nature of this proof, it is possible that more complex techniques could be used to generalize it to a stronger constraint on Euler boxes; this is discussed at the end.
\end{abstract}
\section{Introduction}
The existence of a perfect Euler box (a rectangular cuboid where all sides, face diagonals and space diagonals are integers) is an open question which has been shown to be related to several other open problems \cite{Kanado2023} \cite{Brawer2023}. The problem has been studied since the 1700s, where the eponymous Euler parameterized the construction of imperfect Euler bricks (a perfect Euler box where the space diagonal is allowed to be a noninteger.) Despite the existence of a perfect Euler box still being an open question, many constraints have already been found on perfect Euler boxes, one of which is that the space diagonal of a perfect Euler box cannot be semiprime, i.e. a product of exactly two primes \cite{Korec1992}. This raises the question: is there a similar constraint on the sides of a perfect Euler box? This paper aims to prove the following theorem:
\begin{theorem}
    No Euler boxes exist with a semiprime side length.
\end{theorem}
In Section \ref{variables}, we define the variables and notation that will be used consistently throughout this paper. In Section \ref{constraintsonbandc}, we discuss constraints on the side lengths of an Euler box with semiprime sides. In Sections \ref{case1} and \ref{case2}, we use these constraints to show that no Euler boxes with semiprime sides exist.  In Section \ref{remarks}, we discuss several other related results and approaches which did not succeed. Finally, in Section (to come), we discuss possibilities for generalizing this proof to stronger constraints on Euler boxes. 
\section{Variables and Notation}\label{variables}.
For the duration of this paper, every variable will denote an integer, unless explicitly stated otherwise.
\begin{figure}[H]
  \includegraphics[scale = 0.7]{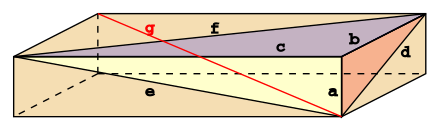}
  \caption{Illustration of an Euler brick. \cite{WikipediaEN:AFM}}
  \label{eulerbrick}
\end{figure}
Suppose we have some rectangular cuboid with positive side lengths $a, b, c$, face diagonals of length $d, e, f$, and space diagonal of length $g$, as shown in Figure \ref{eulerbrick}. If this cuboid is a perfect Euler box, then the following equalities hold:
\begin{align}
    d^2&=a^2+b^2\label{vars1} \\ 
    e^2&=a^2+c^2 \label{vars2} \\
    f^2&=b^2+c^2 \label{vars3}\\
    g^2&=a^2+b^2+c^2 \label{vars4}
\end{align}
\section{Preliminary Constraints on $b$ and $c$}\label{constraintsonbandc}
By equations \ref{vars1} and \ref{vars2}, we have the following equalities:
\begin{align*}
    a^2 &= d^2-b^2\\
    &= (d-b)(d+b)\\
    a^2 &= e^2-c^2\\
    &= (e-c)(e+c)
\end{align*}
We know that $a$ equals $pq$ for some prime $p, q$. Thus, the factors of $a^2$ are given by the following ordered pairs:
\begin{center}
    $(1, p^2q^2),(p,pq^2),(pq,pq),(q,p^2q),(\min{(p^2,q^2)}, \max{(p^2,q^2)}$
\end{center}
where the first element of each ordered pair is less than or equal to the second.\newline 
By the factorization above, we know that the ordered pairs $(d-b,d+b)$ and $(e-c,e+c)$ must be equal to one of the above ordered pairs. (Note that because $b$ and $c$ are positive, all of these ordered pairs have their first element being less than or equal to the second.) In order to narrow down the possibilities for these ordered pairs, we will need to prove a few lemmas.
\begin{lemma}
\label{pairinequality}
The ordered pairs $(d-b,d+b)$ and $(e-c,e+c)$ are not equal.
\end{lemma}\par
\begin{proof}
    Suppose for the sake of contradiction that $(d-b,d+b)$ equals $(e-c,e+c)$. Then we can write:
\begin{align*}
    d+b &= e+c\\
    d-b &= e-c
\end{align*}
Taking the difference of these two equations and simplifying yields:
\begin{align*}
    b &= c
\end{align*}
But by the definition in Section \ref{variables}, we know that the following is true:
\begin{align*}
    f^2 &= b^2+c^2
\end{align*}
Substituting in $b=c$, we get:
\begin{align*}
    f^2 &= 2c^2
\end{align*}
which is a contradiction, as the greatest power of 2 that divides $f^2$ is of the form $2^{2x}$ for some integer x, but the greatest power of 2 that divides $2c^2$ is of the form $2^{2y+1}$ for some integer $y$, and these two quantities can never be equal.
\end{proof} 
\begin{lemma}
\label{pairnotpq}
Neither of the ordered pairs $(d-b,d+b)$ or $(e-c,e+c)$ are equal to $(pq, pq)$.
\end{lemma}\par
\begin{proof}
   Without loss of generality, assume that $(d-b,d+b)$ is equal to $(pq, pq)$. Then $b$ is equal to $0$. But by the definitions in Section \ref{variables}, $b$ is positive, so this is a contradiction. 
\end{proof} 
\begin{lemma}
\label{pairnot1}
Neither of the ordered pairs $(d-b,d+b)$ or $(e-c,e+c)$ are equal to $(1, p^2q^2)$.
\end{lemma}\par
\begin{proof}
    For the sake of contradiction, assume that one of the ordered pairs $(d-b,d+b)$ or $(e-c,e+c)$ are equal to $(1, p^2q^2)$. Without loss of generality, assume that $(d+b,d-b)$ equals $(1, p^2q^2)$. Then we can write:
\begin{align*}
    d+b&=p^2q^2\\
    d-b&=1
\end{align*}
As such, we know that $d$ is equal to $b+1$. By the definitions in Section \ref{variables}, we have:
\begin{align*}
    d^2 &= a^2+b^2\\
    g^2 &= a^2+b^2+c^2
\end{align*}
Substitute in $b+1$ for $d$ to get:
\begin{align*}
    (b+1)^2 &= a^2+b^2
\end{align*}
And substitute $(b+1)^2$ for $a^2+b^2$ to get:
\begin{align*}
    (b+1)^2+c^2&= g^2
\end{align*} 
From equations \ref{vars1} through \ref{vars4}, we have:
\begin{align*}
    b^2+c^2 &= f^2
\end{align*}
Subtract the two equations to yield:
\begin{align*}
    2b+1 &= g^2-f^2\\
    2b+1 &= (f+(g-f))^2-f^2\\
    2b+1 &= 2f(g-f)+(g-f)^2
\end{align*}
Note that $g$ is greater than $f$, so $g-f$ is greater than or equal to $1$. Thus, we have:
\begin{align*}
    2b+1 &\geq 2f+1\\
    b \geq f
\end{align*}
But this is impossible, so we have reached a contradiction. Thus, neither of the ordered pairs $(d-b,d+b)$ or $(e-c,e+c)$ are equal to $(1, p^2q^2)$.
\end{proof}

Using the above lemmas, we can now show that the only possibilities for $(d-b,d+b)$ and $(e-c,e+c)$ are $(p,pq^2)$, $(q,p^2q)$, and $(\min{(p^2,q^2)}, \max{(p^2, q^2)})$. By Lemma \ref{pairinequality}, this leaves us with six possibilities:

\begin{align*}
    (d-b, d+b) &= (p, pq^2) &(e-c, e+c) &= (q, p^2q)\\
    (d-b, d+b) &= (q, p^2q) &(e-c, e+c) &= (p, pq^2)\\
    (d-b, d+b) &= (p, pq^2) &(e-c, e+c) &= (\min{(p^2,q^2)}, \max{(p^2, q^2)})\\
    (d-b, d+b) &= (\min{(p^2,q^2)}, \max{(p^2, q^2)}) &(e-c, e+c)&= (p, pq^2)\\
    (d-b, d+b) &= (q, p^2q) &(e-c, e+c) &= (\min{(p^2,q^2)}, \max{(p^2, q^2)})\\
    (d-b, d+b) &= (\min{(p^2,q^2)}, \max{(p^2, q^2)}) &(e-c, e+c)&= (q, p^2q)
\end{align*}
But because we can interchange $(d-b, d+b)$ and $(e-c, e+c)$ freely, this actually reduces to three equations

\begin{align*}
    (d-b, d+b) &= (p, pq^2) &(e-c, e+c) &= (q, p^2q)\\
    (d-b, d+b) &= (p, pq^2) &(e-c, e+c) &= (\min{(p^2,q^2)}, \max{(p^2, q^2)})\\
    (d-b, d+b) &= (q, p^2q) &(e-c, e+c) &= (\min{(p^2,q^2)}, \max{(p^2, q^2)})\\
\end{align*}

And further, since $p$ and $q$ are interchangeable, we can consolidate this into two cases, with case $1$ being
\begin{align}
    (d-b, d+b) &= (p, pq^2)\label{case1eqstart}\\ 
    (e-c, e+c) &= (q, p^2q)\label{case1eqend}
\end{align}
and case $2$ being 
\begin{align}
    (d-b, d+b) &= (\min{(p^2,q^2)}, \max{(p^2, q^2)})\label{case2eqstart}\\
    (e-c, e+c) &= (p, pq^2)\label{case2eqend}\\
\end{align}
We will first handle case $1$.
\section{Case 1 (p, q symmetric)}\label{case1}
As stated in Equations \ref{case1eqstart} and \ref{case1eqend}, we will assume:
\begin{align*}
    (d-b, d+b) &= (p, pq^2)\\ 
    (e-c, e+c) &= (q, p^2q)
\end{align*}
which, when simplified, yields
\begin{align*}
    b &= \frac{p(q^2-1)}{2}\\
    c &= \frac{q(p^2-1)}{2}
\end{align*}
Using the general properties of Euler boxes as a whole, we can construct an equation relating divisors of $a^2$:
\begin{lemma}
\label{General Case}
Let $d_g = (g+f)$, $d_b = (d+b)$, $d_c = (e+c)$. Then:
\begin{align*}
    \left( \frac{a^2}{d_g} \right) ^2 + 2a^2 + d_g^2 = \left(\frac{a^2}{d_b}\right)^2 + d_b^2 + \left(\frac{a^2}{d_c}\right)^2 + d_c^2
\end{align*}
\end{lemma}
\textit{Proof.} 
Let the variables be as assigned in the statement of Lemma \ref{General Case}. Note that all $d_i$ must divide $a^2$. By the equalities in Section \ref{variables}, we have:
\begin{align*}
    g^2 = \frac{(\frac{a^2}{d_g})^2 + 2a^2 + d_g^2}{4}\\
    b^2 = \frac{(\frac{a^2}{d_b})^2 - 2a^2 + d_b^2}{4}\\
    c^2 = \frac{(\frac{a^2}{d_c})^2 - 2a^2 + d_c^2}{4}\\
\end{align*}
Since $a^2+b^2+c^2 = g^2$, this then resolves to the desired statement. \newline 
By Lemma $\ref{General Case}$, this then yields that for some $d_g \mid a^2$, we have:  
\begin{align*}
    \left(\frac{a^2}{d_g}\right)^2 + 2a^2 + d_g^2 &= p^2q^4 + p^2 + p^4q^2 + q^2\\
    &= (p^2q^2+1)(p^2+q^2)
\end{align*}
The possibilities for $d_g$ are $p^2q^2$, $pq^2$, $pq$, $p^2q$, $p^2$, or $q^2$. $g$ must be strictly greater than $b$ and $c$, so $d_g$ cannot be equal to $d_b$ or $d_c$, and thus cannot be equal to $p^2q$ or $pq^2$. If $d_g = pq$, then $f = 0$, which is impossible, so $d_g \neq pq$. Thus, the only two cases left are $d_g = p^2q^2$, $d_g = p^2$ or $d_g = q^2$. \par 
If $d_g = p^2q^2$, then we have: 
\begin{align*}
    (p^2q^2+1)^2 &= (p^2q^2+1)(p^2+q^2)\\
    p^2q^2+1 &= p^2+q^2\\
    (p^2-1)(q^2-1) &= 0
\end{align*}
which is a contradiction, as $p, q$ are prime. Now without loss of generality, assume $d_g = q^2$. (The equation is the same if $d_g = p^2$.) Then we have: 
\begin{align*}
    (p^2+q^2)^2 &= (p^2q^2+1)(p^2+q^2)\\
    p^2+q^2 &= (p^2q^2+1)\\
    (p^2-1)(q^2-1) &= 0
\end{align*}
which is again a contradiction. Thus, we have solved Case 1.
\section{Case 2 (p, q asymmetric)}\label{case2}
From Equation \ref{case2eqstart} and \ref{case2eqend} we have:
\begin{align*}
    (d-b, d+b) &= (\min{(p^2,q^2)}, \max{(p^2, q^2)})\\
    (e-c, e+c) &= (p, pq^2)\\
\end{align*}
hus, we can write:
$$b = |\frac{p^2-q^2}{2}| = \frac{\pm p^2 \mp q^2}{2}$$
$$c = \frac{p^2q-q}{2}$$
Observe that $a^2 = g^2-f^2$, so $a^2 = (g+f)(g-f)$. Then the pair $(g-f, g+f)$ must be one of the following:
\begin{center}
    $(1, p^2q^2),(p,pq^2),(pq,pq),(q,p^2q),(\min{(p^2,q^2)}, \max{(p^2, q^2)})$
\end{center}
We know that $(g-f, g+f)$ cannot equal either of $(d-b,d+b)$ or $(e-c,e+c)$, because then the diagonal of a right triangle would be equal to one of its legs. Thus, with the assumptions above, this reduces the possibilities to:
\begin{center}
    $(1, p^2q^2),(p,pq^2),(pq,pq)$
\end{center}
Clearly, $(g-f, g+f) \neq (pq, pq)$, because then $f = 0$. Thus, we must only consider the remaining two cases:
\begin{center}
    $(1, p^2q^2),(p,pq^2)$
\end{center}
Suppose that $(g-f, g+f) = (p,pq^2)$. Then we have $g = \frac{pq^2+p}{2}$. Substitute this into the equation $g^2 = a^2+b^2+ c^2$ to yield:
\begin{align*}
    p^2q^4+2p^2q^2+p^2 &= p^4+q^4+p^4q^2+q^2\\
    (pq^2+p)^2 &= (p^4+q^2)(q^2+1)\\
    p^2(q^2+1) &= p^4+q^2\\
    p^4-p^2q^2+q^2-p^2 &= 0\\
    p^2(p^2-1)-q^2(p^2-1) &= 0\\
    (p^2-q^2)(p^2-1) &= 0
\end{align*}
But we have already proven that $p$ and $q$ cannot be equal, and $1$ is not prime, so this is a contradiction. \newline 
Now suppose that $(g-f, g+f) = (1,p^2q^2)$. Then we have $g = \frac{p^2q^2+1}{2}$. Once more, substitute this into $g^2 = a^2+b^2+c^2$ to yield:
\begin{align*}
    p^4q^4+2p^2q^2+1 &= p^4+q^4+p^4q^2+q^2\\
    p^4q^4-p^4-q^4+1 &= p^4q^2+q^2-2p^2q^2\\
    (p^4-1)(q^4-1) &= (q^2)(p^2-1)^2\\
    (p^2+1)(q^4-1) &= (q^2)(p^2-1)\\
    p^2q^4+q^4-p^2-1 &= q^2p^2-q^2\\
    p^2q^4-p^2q^2-p^2+q^4+q^2-1 &= 0\\
    p^2(q^4-q^2-1) + q^4+q^2-1 &= 0
\end{align*}
For all $q > 2$, the left hand side is strictly positive, and it is already known that $q$ is an odd prime, so this is a contradiction. Thus, all possible values for $(d-b,d+b)$ and $(e-c,e+c)$ are invalid, and therefore no Euler boxes with semiprime sides exist.
\section{Remarks} \label{remarks}
Although this bound on its own is relatively limited, the method of assigning the difference of squares to a finite set of factors merits further study. One immediate consequence of this proof is a novel proof of the known result that Euler boxes cannot have an edge of prime length; this can be seen by simply setting $p=1$. Although $1$ is not a prime, the only case where primality was used is in constraining the set of values that each factor in the difference of squares can take. \newline 
This observation can be expressed differently as follows: we can generate possible values for $(d+b, d-b)$ when $a = pq$ by looking at the possible values for $(d+b, d-b)$ when $a=p$ and $a=q$, and pointwise multiplying them together. Explicitly, we can start with
\begin{center}
    $(1, p^2), (p, p)$
\end{center}
and pointwise multiply one element of the former set with one element of the latter set:
\begin{center}
    $(1, q^2), (q, q), (q^2, 1)$
\end{center}
to yield
\begin{center}
    $(1, p^2q^2),(q,p^2q),(p,pq^2), (pq,pq),(\min{(p^2,q^2)}, \max{(p^2, q^2)})$
\end{center}
(with multiplicity). \newline
Reversing this observation yields further, more interesting results. Suppose we were to extend this to considering whether an edge can be a product of exactly three primes, $a = pqr$. What are the possible values that our factors of difference of squares can take? One way to calculate this is by taking any of the values that a factor of difference of squares when $a=pq$ can take:
\begin{center}
    $(1, p^2q^2),(p,pq^2),(pq,pq),(q,p^2q),(\min{(p^2,q^2)}, \max{(p^2, q^2)})$
\end{center}
and pointwise multiplying by any of the values that a factor of difference of squares when $a = r$ can take:
\begin{center}
    $(1, r^2), (r^2, 1), (r, r)$
\end{center}
However, we begin to observe methods of simplification. Observe, for instance, that the pointwise multiplication of $(1, p^2q^2)$ and $(1, r^2)$ yields $(1, p^2q^2r^2)$, but by letting $q' = qr$ we can simplify this to the case of $(d+b, d-b) = (1, p^2q'^2)$, which is within the realm of solutions we have already found. More generally, if $(d+b, d-b) = \Pi_{i = 0}^n (p_i^{a_i})$ where $p_i$ is the $i$th prime factor of $a$ and $a_i \in \{0, 1, 2\}$, then if $a_i = a_j$ are the same, we can set $p_i' = p_ip_j$, reducing the case to a simpler one. This reduces the Euler box problem to a finite set of equations similar to the ones we have solved in this paper; as such, further investigation into these equations could yield stronger constraints. 
\medskip

\printbibliography
\end{document}